\documentclass{amsart}[12pt]
\usepackage{amsmath}
\usepackage{amstext}
\usepackage{amssymb}
\usepackage{amsthm}
\usepackage{amscd}

\usepackage{graphicx}

\usepackage{pstcol}
\usepackage{pstricks}
\usepackage{pst-text}
\usepackage{pst-node}
\usepackage{pst-plot}

\theoremstyle{plain}
\newtheorem{thm}{Theorem}[section]
\newtheorem{lem}[thm]{Lemma}
\newtheorem{prop}[thm]{Proposition}
\newtheorem{cor}[thm]{Corollary}

\theoremstyle{definition}

\theoremstyle{remark}

\DeclareMathOperator{\lcm}{lcm}

\DeclareMathOperator{\Ker}{Ker}
\DeclareMathOperator{\Img}{Im}

\DeclareMathOperator{\Tor}{Tor}
\DeclareMathOperator{\Hred}{\widetilde{H}}
\DeclareMathOperator{\HH}{H}
\DeclareMathOperator{\V}{V}
\DeclareMathOperator{\E}{E}
\DeclareMathOperator{\link}{link}

\begin{document}

\title
[Characteristic-independence of Betti numbers]
{Characteristic-independence of Betti numbers of graph ideals.}
\author{Mordechai Katzman}


\date{\today}

\maketitle

\begin{abstract}
In this paper we study the Betti numbers of Stanley-Reisner ideals
generated in degree 2. We show that the first six Betti numbers do
not depend on the characteristic of the ground field. We also show
that, if the number of variables $n$ is at most 10, all Betti
numbers are independent of the ground field. For $n=11$, there
exists precisely 4 examples in which the Betti numbers depend on
the ground field. This is equivalent to the statement
that the homology of flag complexes with at most 10 vertices is torsion free and that
there exists precisely 4 non-isomorphic flag complexes with 11 vertices whose homology has torsion.

In each of the 4 examples  mentioned above the 8th Betti numbers
depend on the ground field and so we conclude that the highest
Betti number which is always independent of the ground field is
either $6$ or $7$; if the former is true then we show that there
must exist a graph with $12$ vertices whose $7$th Betti number
depends on the ground field.
\end{abstract}

\setcounter{section}{-1}
\section{\bf Introduction}

Throughout this paper $K$ will denote a field.
For any homogeneous ideal $I$ of a polynomial ring $R=K\left[x_1, \dots, x_n\right]$ there exists a \emph{graded} minimal finite free resolution
$$ 0 \rightarrow \bigoplus_j R(- j)^{\beta_{p j}} \rightarrow \dots \rightarrow \bigoplus_j R(-j)^{\beta_{1 j}} \rightarrow
R \rightarrow R/I\rightarrow 0$$
of $R/I$, in which $R(-j)$ denotes the graded free module obtained by shifting the degrees of elements in $R$ by $j$.
The numbers $\beta_{i j}$, which we shall refer to as the $i$th Betti numbers of degree $j$ of $R/I$, are independent of the choice of
graded minimal finite free resolution. We  also define the $i$th Betti number of $I$ as
$\beta_i := \sum \beta_{i j}$.

One of the central problems in Commutative Algebra is the description of minimal resolutions of ideals.
Even when one restricts one's attention to ideals of polynomial rings generated by monomials, the structure of
the resulting resolutions is very poorly understood. There have two main approaches to this problem.
The first is to describe non-minimal free resolutions of these ideals, e.g., the Taylor resolutions
(cf.~\cite{T}) and its generalization, cellular resolutions (cf.~\cite{BS}).
The other approach, which we follow here, has been to
describe the Betti numbers of these minimal resolutions.

It has been known for quite some time that the Betti numbers of monomial ideals may depend
on the characteristic of the ground field
(e.g., see \S 5.4 in \cite{BH1} and section \ref{section4} below.)
The aim of this paper is to investigate this dependence for
Stanley-Reisner rings which are quotients by monomial ideals generated in degree $2$.
In \cite{TH} Naoki Terai and Takayuki Hibi have shown that the third and fourth Betti numbers of these Stanley-Reisner rings
do not depend on the ground field-- this paper extends this result to show that the fifth and sixth
Betti numbers are also independent of the ground field (Theorem \ref{thm2.4} and Corollary \ref{6thBetti}.)
We also show that any such Stanley-Reisner ring whose Betti number depends on the ground field must involve at least
11 variables (Theorem \ref{exampleIsMinimal}) and we list all the minimal examples with 11 variables
(surprisingly, only four such examples exist.)
In these examples the eighth Betti number depends on the ground field and so we conclude that the highest Betti number
which is always independent of the ground field is either $6$ or $7$;
if the former is true then we show that there must exist a
graph with $12$ vertices whose $7$th Betti number depends on the ground field.
Some of the proofs of these results rely on calculations performed by a computer.

Let $G$ be any finite simple graph. We shall always denote the vertex set of $G$ with $\V(G)$ and its edges with $\E(G)$.
Fix a field $K$ and let $K(G)$ be the polynomial ring on the vertices of $G$ over the field $K$.
The graph ideal $I(G)$ associated with $G$ is the ideal of $K(G)$ generated by all degree-2 square-free monomials $u v$
for which $(u,v)\in \E(G)$. It is not hard to see that every ideal in a polynomial ring generated by degree-2 square-free
monomials is of the form $I(G)$ for some graph $G$.

The quotient $K(G)/I(G)$ is a always a Stanley-Reisner ring:
define $\Delta(G)$ to be the simplicial complex on the vertices of $G$ in which
a face consists of a set of vertices, no two joined by an edge. It is easy to see that
$K(G)/I(G)$ coincides with $K[\Delta(G)]$, the Stanley-Reisner ring associated with $\Delta(G)$.
The simplicial complexes of the form $\Delta(G)$ for some graph $G$ are characterised by the fact
that their minimal non-faces have two vertices-- these simplicial complexes are also known as \emph{flag complexes.}

We shall use the following notation and terminology throughout this paper.
For any simple graph $G$, $G^c$ will denote the graph with vertex set $\V(G)$ and edges
$\left\{ (x,y) \,|\, x,y\in\V(G), x\neq y, (x,y)\notin\E(G) \right\}$.

We shall write $\beta_i^K(G)$ and $\beta_{i,d}^K(G)$ for the $i$th Betti number of $K[\Delta(G)]$ and
for the $i$th Betti number of degree $d$ of $K[\Delta(G)]$, respectively.
We may omit the superscript $K$ when the ground field is irrelevant or previously specified.

\section{The Hochster and Eagon-Reiner formulae.}\label{section1}

Recall that for any field $K$ and simplicial complex $\Delta$ the
\emph{Stanley-Reisner ring} $K[\Delta]$ is the quotient of the polynomial ring in the vertices of $\Delta$ with coefficients
in $K$ by the square-free monomial ideal generated by
the product of vertices not in a face of $\Delta$.

The main tool for investigating Betti numbers of a Stanley-Reisner ring $K[\Delta]$ is the following
theorem.
\begin{thm}[Hochster's Formula (Theorem 5.1 in \cite{H})] \label{thm1.0}
The $i$th Betti number of $K[\Delta]$ of degree $d$ is given by
$$\beta_{i,d} = \sum_{W \subseteq V(\Delta), \#W=d} \dim_K \Hred_{d-i-1}(\Delta_W; K)$$
where $V(\Delta)$ is the set of vertices of $\Delta$ and
for any $W \subseteq V(\Delta)$,
$\Delta_W$ denotes the simplicial complex with vertex set $W$ and whose faces are the
faces of $\Delta$ containing only vertices in $W$.

The $i$th Betti number of $K[\Delta]$ is then given by
$$\beta_{i} = \sum_{W \subseteq V(\Delta)} \dim_K \Hred_{\#W-i-1}(\Delta_W; K) .$$
\end{thm}

Notice that when $\Delta=\Delta(G)$ for some graph $G$, we can rewrite the formula above for the Betti numbers as
\begin{equation}\label{HochsterForGraphs}
\beta_{i,d} = \sum_{H \subseteq G\mathrm{\ induced\ } \atop{\#\V(H)=d}} \dim_K \Hred_{d-i-1}(\Delta(H); K) .
\end{equation}

\bigskip
The following is an easy consequence:
\begin{cor}\label{cor1.2}
Let $G$ be any graph.
\begin{enumerate}
\item[(a)] If $H$ is an induced subgraph of $G$ then
$\beta_{i,j}^K(H)\leq \beta_{i,j}^K(G)$ for all fields $K$ and all $i,j\in \mathbb{Z}$.
\item[(b)] $\displaystyle\beta_{i-1,i}^K(G)$ is independent of $K$ and it is non-zero if and only if
$G^c$ contains a disconnected induced subgraph with $i$ vertices.
In particular, the length of the linear strand in a minimal graded resolution of $K[\Delta(G)]$
equals
$$\max \big\{ \V(H)-1 \,|\, H \mathrm{\ is\ a\ disconnected\ induced\ subgraph\ of\ }G^c \big\}.$$
\end{enumerate}
\end{cor}
\begin{proof}
Statement (a) follows immediately from the fact that all summands in (\ref{HochsterForGraphs})
are non-negative.

To prove (b) write
$$\beta_{i-1,i}=\sum_{H \subseteq G\mathrm{\ induced\ } \atop{\#\V(H)=i}} \dim_K \Hred_{0}(\Delta(H); K) $$
and notice that $\Hred_{0}(\Delta(H); K)\neq 0$ if and only if $H^c$ is disconnected, and that
$H^c$ is an induced subgraph of $G^c$ if an only if $H$ is an induced subgraph of $G$.
\end{proof}

\bigskip
The focus of this paper is the study of the dependence of $\beta_{i,j}^K(G)$ on $K$ and we begin
by recording the following basic facts.
\begin{prop}\label{basicFactsOnKDependence}
Let $G$ be any graph.
\begin{enumerate}
\item[(a)] $\beta_{i,j}^K(G)$ depends only on the characteristic of the field $K$.
\item[(b)] $\beta_{i,j}^\mathbb{Q}(G)\leq \beta_{i,j}^{\mathbb{Z}/p\mathbb{Z}}(G)$ for all prime integers $p$.
\item[(c)] $\beta_{i,j}^\mathbb{Q}(G)= \beta_{i,j}^{\mathbb{Z}/p\mathbb{Z}}(G)$ for almost all prime integers $p$.
\item[(d)] $\beta_{i,j}^K(G)$ depends on $K$ if and only if there exists an induced subgraph
$H\subseteq G$ with $j$ vertices
and an $i\geq 1$ for which $\Hred_i(\Delta(H); \mathbb{Z})$ has torsion.
\end{enumerate}
\end{prop}
\begin{proof}

Statement (a) follows from the fact that for any fixed simplicial complex $\Delta$,
$\dim_K \Hred_{i}(\Delta; K)$ depends only on the characteristic of $K$.

Statements (b), (c) and (d) follow from the Universal Coefficient Theorem
(see, for example, Corollary 6.3 in chapter X of \cite{M})
and Hochster's Theorem.
\end{proof}


\bigskip
In \cite{ER} Alexander duality is used to derive a variant of Hochster's Formula.
Recall that for any simplicial complex $\Delta$, the \emph{Alexander Dual} of $\Delta$ is the simplicial complex
defined by
$$\Delta^* := \left\{ F\subseteq \V(\Delta) \,|\, \V(\Delta)-F \notin \Delta \right\} .$$
The link of a face $F\in \Delta$ is defined as the simplicial complex
$$\link_{\Delta} F := \left\{ G\in \Delta \,|\, G\cup F\in \Delta \mathrm{\ and\ } G\cap F=\emptyset \right\} .$$

\begin{thm}[Proposition 1 in \cite{ER}] \label{thm1.1}
The  Betti numbers of $K[\Delta]$ are given by
$$\beta_{i,d} = \sum_{F\in \Delta^*, \ \#(V(\Delta)-F)=d} \dim_K \Hred_{i-2}(\link_{\Delta^*} F; K) $$
and
$$\beta_{i} = \sum_{F\in \Delta^*} \dim_K \Hred_{i-2}(\link_{\Delta^*} F; K) .$$
\end{thm}

When $\Delta=\Delta(G)$ we write $\Delta^*(G)$ for $\left(\Delta(G)\right)^*$.
Notice that faces of $\Delta^*(G)$ are the sets of vertices whose complement
contain two vertices joined by an edge in $G$.
For any $F\in \Delta^*(G)$ the simplicial complex $\link_{\Delta^*} F$ can be easily described as follows:
its maximal faces consist of
$\V(G)-\big(F\cup \{u,v\}\big)$
for all pairs of vertices $u$ and $v$ not in $F$ and which are connected by an edge in $G$.

\bigskip
While Theorem \ref{thm1.1} is essentially identical (via Alexander duality)
to Hochster's Theorem it is often easier to wield as Lemma \ref{LemmaForDegree1}
and its corollary below illustrate.

\bigskip
\begin{lem}\label{LemmaForDegree1}
Let $\Delta$ be a simplicial complex with vertices
$v_1, \dots, v_n$. For any $1\leq i\leq n$ write $\Delta_i$ for the simplex on
$\{v_1, \dots, v_n\} - \{v_i\}$.
Assume that for some $0\leq s\leq n$, $\Delta_{1}, \dots, \Delta_{s}$ are maximal faces of
$\Delta$.
Write $\Delta=\Delta^{(1)} \cup \Delta^{(2)}$
where  $\Delta^{(1)}$ is the sub-complex of $\Delta$ whose maximal faces
are those maximal faces of $\Delta$ which are not among
$\Delta_{1}, \dots, \Delta_{s}$
and where $\Delta^{(2)}= \cup_{i=1}^s \Delta_i$.
If, for some $i\geq 1$, $\dim_K \Hred_{i}(\Delta; K)$ depends on the field $K$,
so does $\displaystyle\dim_K \Hred_{i-s}(\Delta^{(1)}_{\{v_{s+1}, \dots, v_n\}}; K)$.
\end{lem}
\begin{proof}
We proceed by induction on $s$. If $s=0$ the claim is trivial, so assume that $s\geq 1$.
Both $\Delta^{\prime}:=\Delta^{(1)} \cup \Delta_{1} \cup \dots \cup \Delta_{s-1}$
and $\Delta_s$ are acyclic,  the latter because it is a simplex and the former because
$v_s$ is in all its maximal faces and hence is a cone.

The Mayer-Vietoris long exact sequence implies that
$$\Hred_{i}(\Delta; K)\simeq \Hred_{i-1}(\Delta^{\prime} \cap \Delta_s; K)$$
for all $i>1$.
For $i=1$ we obtain the exact sequence
$$0\rightarrow
\Hred_{1}(\Delta; K)\rightarrow \Hred_{0}(\Delta^{\prime} \cap \Delta_s; K) \rightarrow
\Hred_{0}(\Delta^{\prime}; K) \rightarrow \Hred_{0}( \Delta; K) \rightarrow 0 .$$
Since the dimension of the three rightmost $K$-vector spaces is independent of $K$,
$\Hred_{1}(\Delta; K)$ cannot depend on $K$.
We deduce that, if $\dim_K \Hred_{i}(\Delta; K)$ depends on $K$,
$i>1$ and $\dim_K \Hred_{i-1}(\Delta^{\prime} \cap \Delta_s; K)$ also depends on $K$.

We now realise that
$$\Delta^{\prime} \cap \Delta_s =
\big( \Delta^{(1)} \cup  \Delta_1 \cup \dots \cup \Delta_{s-1} \big)_{\big\{v_1, \dots, v_{s-1},v_{s+1}, \dots, v_n\big\}} $$
An application of the induction hypothesis concludes the proof.
\end{proof}

\begin{cor}\label{degree1}
If $G$ contains a vertex $v$ of degree $1$,
then the Betti numbers of $G$ depend on the ground field if and only if those of $G-\{v\}$ do.
\end{cor}
\begin{proof}
If the Betti numbers of $G-\{v\}$ depend on the ground field so do those of $G$ by
Theorem \ref{basicFactsOnKDependence}(d).

Assume now that we can find a counter-example $G$ and pick one with minimal number of vertices.

Let $u$ be the unique neighbour of $v$ in $G$.
Theorem \ref{thm1.1} implies that there exist $i\geq 0$ and  $F\in \Delta^*(G)$ for which
$\dim_K \Hred_{i}(\link_{\Delta^*(G)} F; K)$ depends on $K$.
If $v\in F$,  $\link_{\Delta^*(G)} F=\link_{\Delta^*(G-\{v\})} F-\{v\}$ and the result follows
from the minimality of $G$ together with Theorem \ref{thm1.1}.
If $v\notin F$ but $u\in F$, $v$ is in all maximal faces of $\link_{\Delta^*(G)} F$ and thus the complex is acyclic.

Assume now that $u,v\notin F$, i.e., $u,v \in \link_{\Delta^*(G)} F$.
Notice that $\link_{\Delta^*(G)} F=\link_{\Delta^*(G-F)} \emptyset$ and so the minimality of $G$
implies that $F=\emptyset$.
Write
$\link_{\Delta^*(G)} \emptyset = \Delta^{\prime} \cup \Delta^{\prime\prime}$ where
$\Delta^{\prime}$ is the simplex on the vertices $\V(G)-\{u,v\}$ and
$\Delta^{\prime\prime}$ is the simplicial complex on the vertices $\V(G)$ and whose maximal faces
consist of all $\V(G)- \{x,y\} $ for all edges $(x,y)\in E(G)$ different from $(u,v)$.
Now $\Delta^{\prime}$ and $\Delta^{\prime\prime}$ are acyclic, the former because it is a simplex and the latter because
$v$ is in all its maximal faces and hence is a cone.
The Mayer-Vietoris long exact sequence implies that
$$\Hred_{i}(\link_{\Delta^*(G)} \emptyset; K)\simeq \Hred_{i-1}(\Delta^{\prime} \cap \Delta^{\prime\prime}; K)$$
for all $i>1$. For $i=1$ we obtain the exact sequence
$$0\rightarrow \Hred_{1}(\link_{\Delta^*(G)} \emptyset; K)\rightarrow \Hred_{0}(\Delta^{\prime} \cap \Delta^{\prime\prime}; K)
\rightarrow  \Hred_{0}(\Delta^{\prime\prime}; K) \rightarrow \Hred_{0}(\link_{\Delta^*(G)} \emptyset; K) \rightarrow 0 .$$
Since the dimension of the three rightmost $K$-vector spaces is independent of $K$, so must be the dimension of
$\Hred_{1}(\link_{\Delta^*(G)} \emptyset; K)$.
We deduce that if
$\dim_K \Hred_{i}(\link_{\Delta^*(G)} \emptyset; K)$ depends on $K$ then
$i>1$ and $\dim_K \Hred_{i-1}(\Delta^{\prime} \cap \Delta^{\prime\prime}; K)$ also depends on $K$.

Let $v,u_1, \dots, u_s$ be the neighbours of $u$ among $\V(G)$.
We notice that $\Delta^{\prime} \cap \Delta^{\prime\prime}$ is obtained from
$\Delta^{\prime\prime}$ by removing $u$ and $v$ from all its faces;
so each of the faces $\V(G)- \{u,u_1\} , \dots, \V(G)-\{u,u_s\} $ of
$\Delta^{\prime\prime}$ now correspond to the faces
$$\Delta_1:=\V(G-\{u,v\})-\{u_1\} ,\quad \dots,\quad \Delta_s:=\V(G-\{u,v\})- \{u_s\}
\in \Delta^{\prime} \cap \Delta^{\prime\prime} .$$

We now decompose $\Delta^{\prime}\cap\Delta^{\prime\prime}$
as the union
$\Delta^{(1)} \cup \Delta^{(2)}$
where $\Delta^{(2)}=\Delta_1 \cup \dots \cup \Delta_s$ and
$\Delta^{(1)}$ is the sub-simplicial complex of $\Delta^{\prime\prime}_{\V(G)-\{u,v\}}$
whose maximal faces are those maximal faces of $\Delta^{\prime\prime}_{\V(G)-\{u,v\}}$
which are not among $\Delta_1, \dots , \Delta_s$.
Now Lemma \ref{LemmaForDegree1} implies that
$\displaystyle\dim_K \Hred_{i-s-1}(\Delta^{(1)}_{\big\{\V(G)-\{v,u,u_1, \dots, u_s\} \big\}}; K)$
depends on $K$.
But it is not hard to see that
$\Delta^{(1)}_{\big\{\V(G)-\{v,u,u_1, \dots, u_s\} \big\}}=\link_{\Delta^*(G-\{v\})} \big\{u,u_1, \dots, u_s\big\}$,
and so we are done by Theorem \ref{thm1.1}.
\end{proof}

\bigskip
In what follows we shall also need the following theorem proved in \cite{JK} and in \cite{J}.
\begin{thm}[\cite{JK} and \cite{J}]\hfil\nobreak\label{DegreeAtMost2}
\begin{enumerate}
\item[(a)] Let $G_1$ and $G_2$ be disjoint graphs and let $G=G_1 \cup G_2$.
The Betti numbers $\beta_{i,j}^K(G)$ are independent of $K$ if and only if
$\beta_{i,j}^K(G_1)$ and $\beta_{i,j}^K(G_2)$ are independent of $K$.

\item[(b)] If the vertices of $G$ have degree at most 2 then the Betti numbers of $K[\Delta(G)]$ do not depend on $K$.
Consequently, $\Hred_j\left(\Delta(G); \mathbb{Z}\right)$ and $\Hred_j\left(\Delta^*(G); \mathbb{Z}\right)$
are torsion free for all $j\in \mathbb{Z}$.
\end{enumerate}
\end{thm}

\bigskip
One of the aims of the study of the Betti numbers of graph ideals is the search
for their combinatorial significance.  Corollary \ref{cor1.2} is an example of such an interpretation
(see \cite{J} and \cite{JK} for more results of this type.)

One could think that, if these Betti numbers can be interpreted purely in terms of the combinatorial structure of $G$,
the choice of ground field $K$ should not affect the values of the Betti numbers.
This is not the case, as we shall see in section \ref{section4}.

\section{Applications of Taylor's resolution.}\label{section2}

Let $K$ be a field, $m_1, \dots, m_n$ any monomials in $R=K[x_1, \dots, x_s]$ and let $I$ be the ideal
generated by $m_1, \dots, m_n$.
In \cite{T} Diana Taylor produced an explicit construction of a free (but seldom minimal) resolution
for $R/I$ which we now describe.

For every $0\leq i\leq n$ define $\mathcal{G}_i$ to be the set of length-$i$ subsequences
$(m_{j_1} , m_{j_2}, \dots, m_{j_i})$ of $(m_1, \dots, m_n)$.

For every $1\leq i\leq n$ let $T_i$ be the free $R$-module whose free generating set is $\mathcal{G}_i$
and define $T_0=R$.
Now for all $i\geq 1$ define $\partial_i : T_i \rightarrow T_{i-1}$ by specifying
\begin{equation}\label{TaylorDifferential}
\partial_i(m_{j_1}, \dots, m_{j_i})=
\sum_{k=1}^i (-1)^k
\frac
 {\lcm \{m_{j_1}, \dots, m_{j_i}\}}
 {\lcm \{m_{j_1}, \dots, m_{j_{k-1}},m_{j_{k+1}},\dots m_{j_i}\}} (m_{j_1}, \dots, m_{j_{k-1}},m_{j_{k+1}},\dots m_{j_i}) .
\end{equation}

If we further declare the degree of each free generator $g\in \mathcal{G}_i$ to be
$\deg \lcm g$, $T_\bullet$ becomes a \emph{graded} free resolution.

Although $T_\bullet$ is not minimal, we may use it to compute the $i$th Betti numbers of degree $d$ of $R/I$ as
$$\Tor^R_i\big(R/I,R/(R x_1 + \dots R x_s) \big)_d = \HH_i \big(T_\bullet \otimes_R R/(R x_1 + \dots R x_s)\big)_d .$$

The following is an easy observation following from this construction:
\begin{prop}\label{zeroBeyondDegLcm}
Let
$D=\max_{g\in \mathcal{G}_i} \deg \lcm g$.
The $i$th Betti number of degree $d$ of $R/I$ vanishes for all $d>D$.
\end{prop}

\bigskip
Now we restrict our attention to Taylor resolutions of graph ideals.
Fix an ordering of the edges of $G$, $e_1, \dots, e_E$.
We can think now of
$T_i$ as being the free $R$-module whose free generators consist of sequences
$(e_{j_1}, \dots, e_{j_i})$
of $i$ edges in $G$ where $j_1< \dots < j_i$
and we can rewrite (\ref{TaylorDifferential}) as
\begin{equation}\label{TaylorDifferentialG1}
\partial_i (e_{j_1}, \dots, e_{j_i})=
\sum_{k} (-1)^k  \mu_k (e_{j_1}, \dots, e_{j_{k-1}}, e_{j_{k+1}}, \dots e_{j_i}) .
\end{equation}
where $\mu_k$ is the product of the vertices in $e_{j_k}$ which are not in any of
$e_{j_1}, \dots, e_{j_{k-1}},e_{j_{k+1}}, \dots, e_{j_i}$.

Let $J$ be the ideal of $R(G)$, the polynomial ring over $K$ in the vertices of $G$,
generated by the vertices of $G$.

Notice that, after tensoring with $R(G)/J$,
$\partial_i (e_{j_1}, \dots, e_{j_i})$ vanishes unless there exists a $1\leq k\leq i$
such that both vertices in $e_{j_k}$ occur in $e_{j_1}, \dots, e_{j_{k-1}},e_{j_{k+1}}, \dots, e_{j_i}$.
So the differentials in $T_\bullet \otimes_R R(G)/J$ are defined by
\begin{equation}\label{TaylorDifferentialG}
\overline{\partial}_i \{e_1, \dots, e_i\}=
\sum_{\mathrm{vertices\ of\ }e_{j_k}\mathrm{\ are\ in \ }\atop{e_{j_1}, \dots, e_{j_{k-1}},e_{j_{k+1}}, \dots, e_{j_i}}}
(-1)^k  (e_{j_1}, \dots, e_{j_{k-1}}, e_{j_{k+1}}, \dots e_{j_i})  .
\end{equation}

\begin{lem}\label{BettiZeroBeyond2i}
For any graph $G$ and any $i\geq 1$, $\beta_{i,d}(G)=0$ for all $d>2i$ and
$\beta_{i,2i}$ is the number of induced subgraphs of $G$ consisting of $i$ disjoint edges.
\end{lem}
\begin{proof}
The first statement is a consequence of Proposition \ref{zeroBeyondDegLcm}.

Notice that the degree-$2i$ free generators of $T_{i}$ are
those sets of $i$ edges which together contain $2i$ vertices, and that the only such
sets of edges are sets of $i$ disjoint edges.
An easy examination of (\ref{TaylorDifferentialG}) shows that, for such free generator $g$,
the image of $\partial_i g$ in $T_i \otimes_R R(G)/J$ vanishes.
Also, if these $i$ disjoint edges $\{e_{j_1}, \dots, e_{j_i}\}$
do not form an induced subgraph of $G$, i.e., if there exists another edge
$e$ whose both vertices occur in $\{e_{j_1}, \dots, e_{j_i}\}$ then,
working modulo $J$,
$\overline{\partial_{i+1}} \{e_{j_1}, \dots, e_{j_i},e\}=(e_{j_1}, \dots, e_{j_i})$.
Finally, if the $i$ disjoint edges $\{e_{j_1}, \dots, e_{j_i}\}$ form an induced subgraph,
the generator $(e_{j_1}, \dots, e_{j_i})$ cannot occur in image of
$\overline{\partial_{i+1}} (t)$ for any $t\in T_{i+1}$.
To see this note that it can only occur in $\overline{\partial_{i+1}} \{e_{j_1}, \dots, e_{j_i},e\}$
for some edge $e$ and that the fact that edges $\{e_{j_1}, \dots, e_{j_i}\}$ form an induced subgraph of $G$
implies that at least one of the vertices in $e$ does not occur in $\{e_{j_1}, \dots, e_{j_i}\}$ and, therefore,
the coefficient of $(e_{j_1}, \dots, e_{j_i})$ in $\overline{\partial_{i+1}} \{e_1, \dots, e_i,e\}$ is zero.
We now conclude that $\HH_i(T_\bullet \otimes_R R(G)/J)$ has a $K$-basis consisting of all induced subgraphs of $G$
consisting of $i$ disjoint edges.
\end{proof}

\begin{lem}[see also Lemma 2.1 in \cite{TH}]\label{lem2.0}
Let $G$ be a graph with $n$ vertices. If $n<2(j+1)$ then
$\Hred_{j}\left(\Delta(G); \mathbb{Z}\right)=0$ and, if
$n=2(j+1)$, $\Hred_{j}\left(\Delta(G); \mathbb{Z}\right)=0$ unless
$G$ consists of $j+1$ disjoint edges.
\end{lem}
\begin{proof}
To prove the first statement rewrite $n<2(j+1)$ as $n>2(n-j-1)$ and notice that Lemma \ref{BettiZeroBeyond2i}
implies that $\beta_{n-j-1,n}(G)=0$ and that
Hochster's Theorem shows that for any field $K$
$$0=\beta_{n-j-1,n}^K(G)=\Hred_{j}(\Delta(G); K).$$

If $n=2(j+1)$, for any field $K$
$$\beta_{n-j-1,n}^K(G)=\beta_{j+1,2(j+1)}^K(G)=\Hred_{j}(\Delta(G); K) $$
and the result follows from the second statement in Lemma \ref{BettiZeroBeyond2i}.
\end{proof}

\bigskip
We shall also need the following result
\begin{lem}\label{Betti i,2i-1}
For any graph $G$ and any $i\geq 1$, $\beta^K_{i,2i-1}(G)$ does not depend on $K$.
\end{lem}
\begin{proof}
If a counter-example exists, Proposition \ref{basicFactsOnKDependence}(d) implies that
we may, and shall, choose the counter-example $G$ to have $2i-1$ vertices.
Pick such a counterexample with minimal $i$.
Theorem 2.4 in \cite{TH} (and also Theorem \ref{thm2.4} in this paper) implies that $i\geq 4$.

Pick a free generator in $T_i \otimes_R R(G)/J$ consisting of a
set of $i$ edges involving $2i-1$ vertices, i.e., a subgraph $G^\prime$ of $G$ of the form
$$
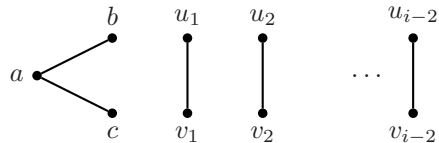
\psset{xunit=1}
\psset{yunit=0.5}
\begin{pspicture}(-1,-0.5)(6,2.5)
\dotnode(0,1){a} \uput[180](0,1){$a$}
\dotnode(1,0){b} \uput[-90](1,0){$c$}
\dotnode(1,2){c} \uput[90](1,2){$b$}
\dotnode(2,2){u1}\uput[90](2,2){$u_1$}
\dotnode(2,0){v1}\uput[-90](2,0){$v_1$}
\dotnode(3,2){u2}\uput[90](3,2){$u_2$}
\dotnode(3,0){v2}\uput[-90](3,0){$v_2$}
\uput[0](4,1){$\dots$}
\dotnode(5,2){uk}\uput[90](5,2){$u_{i-2}$}
\dotnode(5,0){vk}\uput[-90](5,0){$v_{i-2}$}

\ncline{-}{a}{b}
\ncline{-}{a}{c}
\ncline{-}{u1}{v1}
\ncline{-}{u2}{v2}
\ncline{-}{uk}{vk}
\end{pspicture}
$$
and pick this generator so that its image in $\Ker \overline{\partial}_{i}/\Img \overline{\partial}_{i+1}$ is not zero.

As in the proof of Lemma \ref{BettiZeroBeyond2i}, the only edges in $G$ among the vertices
$u_1, \dots, u_{i-2}, v_1, \dots, v_{i-2}$ are $(u_1,v_1), \dots, (u_{i-2},v_{i-2})$.
Also, none of these vertices is joined by an edge to $a$, otherwise, if, say, $(u_1,a)\in \E(G)$, then
$$\overline{\partial_{i+1}} \big\{ (a,b),(a,c),(a,u_1),(u_1,v_1), \dots, (u_{i-2},v_{i-2}) \big\}=
\pm \big\{ (a,b),(a,c),(u_1,v_1), \dots, (u_{i-2},v_{i-2}) \big\}$$
contradicting the fact that the image of $G^\prime$ in $\Ker \overline{\partial}_{i}/\Img \overline{\partial}_{i+1}$ is not zero.

We now proceed by examining an exhaustive set of cases.

\noindent\emph{\underline{Case I}: For some $1\leq j\leq i-2$ both vertices $u_j$ and $v_j$ have degree 1.}
Assume with no loss of generality that $j=1$.
Hochster's formula gives
$$\beta^K_{i,2i-1}(G) =  \dim_K \Hred_{i-2}(\Delta(G); K) .$$
But $\Delta(G)$ is the suspension of $\Delta(G-\{u_1,v_1\})$ and hence
$\Hred_{i-2}(\Delta(G); K)\cong \Hred_{i-3}(\Delta(G-\{u_1,v_1\}); K)$.
Another application of Hochster's formula gives
$$\beta^K_{i-1,2(i-1)-1}(G-\{u_1,v_1\})=\dim_K \Hred_{i-3}(\Delta(G-\{u_1,v_1\}); K)$$
which, by the minimality of $G$, is independent of $K$. We deduce that $\beta^K_{i,2i-1}(G)$ is also independent of $K$.

\noindent\emph{\underline{Case II}: there exist $1\leq j_1, j_2 \leq i-2$, $j_1\neq j_2$ for which
$$\{(u_{j_1},c),(v_{j_1},c)\} \cap \E(G) \neq \emptyset \text{ and } \{(u_{j_2},b),(v_{j_2},b)\} \cap \E(G) \neq \emptyset .$$
}

Assume with no loss of generality that $j_1=1$, $j_2=2$ and that $(v_1,c),(u_2,b) \in \E(G)$. Now
\begin{eqnarray*}
\overline{\partial_{i+1}}
\lefteqn{\big\{ (a,b),(a,c),(v_1,c), (u_1,v_1), \dots, (u_{i-2},v_{i-2}) \big\}}\\
&=&\pm \big\{ (a,b),(a,c), (u_1,v_1), \dots, (u_{i-2},v_{i-2}) \big\}
\pm \big\{ (a,b),(v_1,c), (u_1,v_1), \dots, (u_{i-2},v_{i-2}) \big\}
\end{eqnarray*}
But the edges $(u_2,v_2)$ and $(a,b)$ are joined by $(u_2,b)$
and so
$$\overline{\partial_{i+1}}
\big\{ (a,b),(v_1,c), (u_1,v_1), (u_2,b) \dots, (u_{i-2},v_{i-2}) \big\}=
\big\{ (a,b),(v_1,c), (u_1,v_1), \dots, (u_{i-2},v_{i-2}) \big\}$$
and the image of
$$ \big\{ (a,b),(v_1,c), (u_1,v_1), \dots, (u_{i-2},v_{i-2}) \big\}$$
in $\Ker \overline{\partial}_{i}/\Img \overline{\partial}_{i+1}$ is zero and so the image of $G^\prime$ in $\Ker \overline{\partial}_{i}/\Img \overline{\partial}_{i+1}$ vanishes,
a contradiction.

Since $G^\prime$ contains at least two isolated edges, if none of the two cases above hold,
at least one of the vertices $b$ or $c$ must not be a neighbour of any of $u_1,v_1, \dots, u_{i-2},v_{i-2}$.

\noindent\emph{\underline{Case III}: $\deg b=1$ or $\deg c=1$.}
With no loss of generality assume that the former occurs.
An application of Theorem \ref{thm1.1} gives
$$\beta^K_{i,2i-1}(G) =
\dim_K \Hred_{i-2}(\link_{\Delta^*(G)} \emptyset; K)=
\dim_K \Hred_{i-2}(\Delta^*(G); K).$$

Let $\Delta_1$ be the simplex with vertices $\V(G)-\{a,b\}$ and let
$\Delta_2$ be the simplicial complex with vertex-set $\V(G)$ and whose maximal faces are
$$\big\{ V(G) -\{x,y\} \,|\, (x,y)\in \E(G)-\{(a,b)\} \big\} .$$
Notice that $\Delta_1$ and $\Delta_2$ are acyclic, the latter because $b$ is in all maximal faces.
It follows from the discussion after Theorem \ref{thm1.1} that
$\Delta^*(G)=\Delta_1 \cup \Delta_2$ and the corresponding Mayer-Vietoris long exact sequence gives
$$\Hred_{i-2}(\Delta^*(G); K)\cong \Hred_{i-3}(\Delta_1 \cap \Delta_2; K). $$
But $\Delta_1 \cap \Delta_2$ is a simplicial complex with vertex-set $X:=\V(G)-\{a,b\}$ and whose maximal faces are
$$ \big\{X-\{c\}\big\} \cup \big\{ X-\{u_j,v_j\} \,|\, 1\leq j\leq i-2 \big\} .$$
Let $\Delta_3$ be the simplex with vertices $X-\{c\}$ and let
$\Delta_4$ be the simplicial complex with vertex-set $X$ and whose maximal faces are
$$\big\{ X-\{u_j,v_j\} \,|\, 1\leq j\leq i-2 \big\} ;$$
we now can write $\Delta_1 \cap \Delta_2=\Delta_3 \cup \Delta_4$.
We apply Lemma \ref{LemmaForDegree1} to deduce that, if
$\dim_K \Hred_{i-4}((\Delta_4)_{X-\{c\}}; K)$ is independent of $K$, so is
$\dim_K \Hred_{i-3}(\Delta_3 \cup \Delta_4; K)=\dim_K \Hred_{i-3}(\Delta_1 \cap \Delta_2; K)$.

Let $H$ be the induced subgraph of $G$ with vertex set $X-\{c\}$, i.e.,
the disjoint union of the edges $\{u_1,v_1\}, \dots, \{u_{i-2},v_{i-2}\}$
and notice that $\Delta^*(H)=(\Delta_4)_{X-\{c\}}$.

Alexander Duality (cf. Theorem 6.2 in \cite{B}) implies that
$\Hred_{i-4}(\Delta^*(H); K)\cong \Hred^{i-3}(\Delta(H); K)$,
but $\Delta(H)$ is a sphere (it is a repeated suspension of a 0-sphere,)
and so  its cohomology is independent of $K$.

\end{proof}

Another consequence of equation (\ref{TaylorDifferentialG1}) is the following.

\begin{prop}\label{StarsAndProjectiveDimension}
Assume that the graph $G$ contains an induced subgraph $H$
with $i$ edges and $d$ vertices in which all edges contain a vertex of degree one.
Then $\beta_{i,d}^K(G)\neq 0$ for all fields $K$.
\end{prop}
\begin{proof}
In view of Proposition \ref{cor1.2}(a) it is enough to show that
$\beta_{i,d}^K(H)\neq 0$.

Let $e_1, \dots, e_i$ be all the edges of $H$, let $T_\bullet$ be the Taylor resolution of $K[\Delta(H)]$
and consider the free generator  $(e_1, \dots, e_i)$ in  $T_i$.
The degree of the generator is $d$,
$\overline{\partial}_i(e_1, \dots, e_i)=0$ and $\overline{\partial}_{i+1}=0$ so
 $(e_1, \dots, e_i)$ represents a non-zero element in $\Ker \overline{\partial}_{i}/\Img \overline{\partial}_{i+1}$
and, therefore, $\beta_{i,d}^K(H)\neq 0$.
\end{proof}

\section{Low Betti numbers of graph ideals.}\label{section3}

In \cite{TH} it is shown that the third and fourth betti numbers of $K[\Delta(G)]$ do not depend on $K$.
The main result in this section, Theorem \ref{thm2.4}, extends this result and shows that the fifth Betti number
of $K[\Delta(G)]$ does not depend on $K$ either.
We shall see later that the sixth Betti number also does not depend on $K$

\begin{lem}\label{DegreeAtMost3}
Let $G$ be a graph with $n$ vertices.
If the vertices of $G$ have degree at most 3 then $\Hred_{j}\left(\Delta(G); \mathbb{Z}\right)$ is torsion-free for all
$j\geq n-6$.
\end{lem}
\begin{proof}
Let $G$ be a counterexample with minimal number of vertices $n$.
If  all vertices in $G$ have degree at most 2, the result is a consequence of Theorem \ref{DegreeAtMost2}(b).
Assume that we can find a vertex $v$ in $G$ whose degree is $3$
and let $\left\{v_1,\dots, v_{n-4}\right\}$ be the set of vertices in $G$ which are not neighbours of $v$.
Let $H$ be the induced subgraph of $G$ with vertices $\left\{v_1,\dots, v_{n-4}\right\}$ and
let $H^\prime$ be the induced subgraph of $G$ with vertices $\left\{v,v_1,\dots, v_{n-4}\right\}$.

Notice that $\Delta(G)=\Delta(G-\{v\}) \cup \Delta(H^\prime)$ and that $\Delta(H^\prime)$ is a cone and hence acyclic.
Consider the following Mayer-Vietoris exact sequence
\begin{eqnarray}\label{MV1}
\dots \rightarrow \Hred_{j}\Big(\Delta(G-\{v\}) \cap \Delta(H^\prime); \mathbb{Z}\Big) \rightarrow
\Hred_{j}\Big(\Delta(G-\{v\}); \mathbb{Z}\Big) \rightarrow\\
\nonumber \Hred_{j}\Big(\Delta(G); \mathbb{Z}\Big) \rightarrow
\Hred_{j-1}\Big(\Delta(G-\{v\}) \cap \Delta(H^\prime); \mathbb{Z}\Big) \rightarrow \dots
\end{eqnarray}
To show that $\Hred_{j}\left(\Delta(G); \mathbb{Z}\right)$ is torsion-free for all $j\geq n-6$
it is enough to show that for all $j\geq \max\{1,n-6\}$,
$\Hred_{j}\big(\Delta(G-\{v\}) \cap \Delta(H^\prime); \mathbb{Z}\big)=0$ and that
$\Hred_{j}\big(\Delta(G-\{v\}); \mathbb{Z}\big)$ and
$\Hred_{j-1}\big(\Delta(G-\{v\}) \cap \Delta(H^\prime); \mathbb{Z}\big)$ are torsion-free.

Notice that $\Delta(G-\{v\}) \cap \Delta(H^\prime)=\Delta(H)$ and so the minimality of $n$ implies that
$\Hred_{j}\big(\Delta(H^\prime); \mathbb{Z}\big)$ and
$\Hred_{j-1}\big(\Delta(G-\{v\}) \cap \Delta(H^\prime); \mathbb{Z}\big)$ are torsion-free
for all $j\geq n-6$.

Whenever $n\geq 7$ and $j\geq n-6$ we have $2(j+1)\geq 2n-10>n-4$ and, if
we apply Lemma \ref{lem2.0} to $\Delta(G-\{v\}) \cap \Delta(H^\prime)=\Delta(H)$,
we see that whenever $n\geq 7$ and $j\geq n-6$ we have
$\Hred_{j}\left(\Delta(H); \mathbb{Z}\right)=0$.
On the other hand, if $n<7$,
$H$ contains at most two vertices and clearly $\Hred_{j}\big(\Delta(H); \mathbb{Z}\big)=0$ for all $j>0$.

\end{proof}

\begin{lem}\label{DegreeAtLeastn-4}
\begin{enumerate}
\item[(a)] Assume that $G$ is a graph with $n$ vertices which contains a vertex $v$ of degree $n-1$.
Let $\beta_i$ and  $\beta^\prime_i$ be the $i$-th Betti numbers of
$K\left[\Delta(G)\right]$ and $K\left[\Delta(G-\{v\})\right]$, respectively.
Then for all $i>1$
$$\beta_i = \beta_i^\prime + \beta^\prime_{i-1} + {{n-1} \choose{i}} .$$

\item[(b)] Assume that $G$ is a graph with $n$ vertices which contains a vertex $v$ of degree at least $n-4$.
The Betti numbers of
$K\left[\Delta(G)\right]$  are independent of the characteristic of $K$
if and only if
the Betti numbers of $K\left[\Delta(G-\{v\})\right]$ are independent of the characteristic of $K$.
\end{enumerate}

\end{lem}
\begin{proof}
Let $v$ be a vertex of $G$ of degree $n-1$.
We use Hochster's formula for the Betti numbers of $K\left[\Delta(G)\right]$ as follows
\begin{eqnarray*}
\beta_i & = & \sum_{V \subseteq \V(G)} \dim_K \Hred_{\#V-i-1}(\Delta(G)_V; K)\\
        & = & \sum_{V \subseteq \V(G), v\notin V} \dim_K \Hred_{\#V-i-1}(\Delta(G)_V; K) + \sum_{V \subseteq \V(G), v\in V} \dim_K \Hred_{\#V-i-1}(\Delta(G)_V; K)\\
        & = & \beta_i^\prime + \sum_{V \subseteq \V(G), v\in V} \dim_K \Hred_{\#V-i-1}(\Delta(G)_V; K)
\end{eqnarray*}
Notice that the only face of $\Delta(G)$ which contains $v$ is the $0$-dimensional face $\{ v \}$.
So, if $V \subseteq \V(G)$ and $v\in V$, for all $i>1$
$$
\dim_K \Hred_{\#V-i-1}(\Delta(G)_V; K) =
\left\{
\begin{array}{lll}
1+ \dim_K \Hred_{\#V-i-1}(\Delta(G)_{V-\{v\}}; K), & \mathrm{if\ } \#V-i-1=0\\
\dim_K \Hred_{\#V-i-1}(\Delta(G)_{V-\{v\}}; K), & \mathrm{otherwise.}\\
\end{array}
\right.
$$
and so
$$
\begin{array}{l}
\displaystyle\sum_{V \subseteq \V(G)\atop{v\in V}}\dim_K \Hred_{\#V-i-1}\big(\Delta(G)_V; K\big) \\
\displaystyle=\sum_{V \subseteq \V(G)\atop{v\in V, \#V=i+1}} 1+\dim_K \Hred_{\#V-i-1}\big(\Delta(G)_{V}; K\big) +
\sum_{V \subseteq \V(G)\atop{v\in V, \#V\neq i+1}} \dim_K \Hred_{\#V-i-1}\big(\Delta(G)_{V}; K\big)\\
\displaystyle=\sum_{U \subseteq \V(G-\{v\}) \atop{\#U=i}} 1+\dim_K \Hred_{\#U-i}\big(\Delta(G-\{v\})_{U}; K\big) +
\sum_{U \subseteq \V(G-\{v\}) \atop{\#U\neq i}} \dim_K \Hred_{\#U-i}\big(\Delta(G-\{v\})_{U}; K\big)\\
\displaystyle= {{n-1} \choose{i}} + \beta^\prime_{i-1,i} + \sum_{j\neq i} \beta^\prime_{i-1,j}
= {{n-1} \choose{i}} + \beta^\prime_{i-1} .
\end{array}
$$
We now obtain for all $i>1$
$$\beta_i = \beta_i^\prime + \beta^\prime_{i-1} + {{n-1} \choose{i}} .$$

It is enough to show that if the Betti numbers of $K\left[\Delta(G-\{v\})\right]$ are independent of the characteristic of $K$
so are those of $K\left[\Delta(G)\right]$.
Pick a counter-example $G$ with minimal number of vertices $n$.
When $v$ has degree $n-1$ (b) follows easily from (a).
Assume now that $v$ has degree at most $n-2$.
By Theorem \ref{DegreeAtMost2} we may assume that $G$ is connected.

Proposition \ref{basicFactsOnKDependence} implies that the Betti numbers of
$K[\Delta(G)]$ are independent of $K$ if and only if the $\mathbb{Z}$-module
$\HH_i(\Delta(G); \mathbb{Z})$ has no torsion for all $i\geq 1$.

Let $v_1, \dots, v_s\in\V(G)$ be the non-neighbours of $v$ and let $H$ be the subgraph induced by them.
Since $s\leq 3$, $\Hred_i(\Delta(H); \mathbb{Z})=0$ for all $i>0$;

Let $H^\prime$ be the induced subgraph of $G$ with vertices $v,v_1, \dots, v_s$.
Notice that $\Delta(H^\prime)$ is a cone and hence acyclic.
We have $\Delta(G)=\Delta(G-\{v\}) \cup \Delta(H^\prime)$ and
$\Delta(G-\{v\}) \cap \Delta(H^\prime)=\Delta(H)$.
The corresponding Mayer-Vietoris exact sequence gives an isomorphism
$
\Hred_i(\Delta(G-\{v\}); \mathbb{Z})  \cong
\Hred_i(\Delta(G); \mathbb{Z})
$
for all $i>1$.
We also obtain the following Mayer-Vietoris exact sequence
$$
0 \rightarrow
\Hred_1(\Delta(G-\{v\}); \mathbb{Z})  \rightarrow
\Hred_1(\Delta(G); \mathbb{Z}) \rightarrow
\Hred_0(\Delta(H); \mathbb{Z}) \rightarrow
0 .
$$
Since both
$\Hred_1(\Delta(G-\{v\}); \mathbb{Z})$ and
$\Hred_0(\Delta(H); \mathbb{Z})$ are torsion-free,
so is $\Hred_1(\Delta(G); \mathbb{Z})$.
\end{proof}

\begin{lem}\label{lem2.4}
Let $G$ be a graph which contains a vertex of degree $\delta\geq 4$.
If
$\Hred_{j-1} \left(\Delta^*(G-\{v\}); \mathbb{Z}\right)$ is torsion-free for all $j\leq 3$ then
$\Hred_{j} \left(\Delta^*(G); \mathbb{Z}\right)$ is torsion-free for all $j\leq 3$.
\end{lem}
\begin{proof}
Let $v$ be a vertex of $G$ of degree $\delta$ and
let $v_1,\dots ,v_\delta$ be the neighbours of $v$.
If $\left\{v, v_1,\dots ,v_\delta\right\} = \V(G)$, we are done by Lemma \ref{DegreeAtLeastn-4},
so we may assume that there exists a vertex $w\in \V(G)-\left\{v, v_1,\dots ,v_\delta\right\}$.

Let $G_1$ and $G_2$ be subgraphs of $G$ which contain all its vertices;
let the edges of $G_1$ be
\goodbreak
$\left\{(v, v_1), \dots, (v,v_\delta)\right\}$
and let the edges of $G_2$ be $E(G)-E(G_1)$.

It is not hard to see that  $\Delta^*(G)=\Delta^*(G_1) \cup \Delta^*(G_2)$.
Consider the following Mayer-Vietoris exact sequence
\begin{eqnarray*}\label{MV2}
\dots \rightarrow
\Hred_{j}\left(\Delta^*(G_1); \mathbb{Z}\right) \oplus \Hred_{j}\left(\Delta^*(G_2); \mathbb{Z}\right) \rightarrow
\Hred_{j}\left(\Delta^*(G_1) \cup \Delta^*(G_2); \mathbb{Z}\right) \rightarrow \\
\Hred_{j-1}\left(\Delta^*(G_1) \cap \Delta^*(G_2); \mathbb{Z}\right) \rightarrow
\Hred_{j-1}\left(\Delta^*(G_1); \mathbb{Z}\right) \oplus \Hred_{j-1}\left(\Delta^*(G_2); \mathbb{Z}\right) \rightarrow \dots
\end{eqnarray*}
Since $w$ is in all maximal faces of $\Delta^*(G_1)$ we have
$\Hred_{j}\left(\Delta^*(G_1); \mathbb{Z}\right)=0$ for all $j$ and
since $v$ is in all maximal faces of $\Delta^*(G_2)$ we have
$\Hred_{j}\left(\Delta^*(G_2); \mathbb{Z}\right)=0$ for all $j$.
So, for all $j$,
\begin{equation}\label{eqn1}
\Hred_{j}\left(\Delta^*(G_1) \cup \Delta^*(G_2); \mathbb{Z}\right) \cong
\Hred_{j-1}\left(\Delta^*(G_1) \cap \Delta^*(G_2); \mathbb{Z}\right) .
\end{equation}

Notice that $\Delta^*(G_1),\Delta^*(G_2)\subseteq \Delta^*(G)$ and so $\Delta^*(G_1) \cap \Delta^*(G_2)\subseteq \Delta^*(G)$.
Furthermore, $v$ is not in any maximal face of $\Delta^*(G_1) \cap \Delta^*(G_2)$, so
$\Delta^*(G_1) \cap \Delta^*(G_2) \subseteq \Delta^*(G-\{v\})$.
We now show that for all $d\leq \delta-2$, any $d$-dimensional face $f$ of $\Delta^*(G-\{v\})$ is also a face in
$\Delta^*(G_1) \cap \Delta^*(G_2)$.
Any such face $f$ must exclude two vertices of $G-\{v\}$ joined by an edge in $G-\{v\}$, so $f\in \Delta^*(G_2)$.
Also, as $\dim f \leq \delta-2$, $f$ has to exclude at least one $v_1, \dots, v_\delta$, and, therefore,
$f\in \Delta^*(G_1)$.
We may now deduce that
\begin{equation}\label{eqn2}
\Hred_{j-1}\left(\Delta^*(G_1) \cap \Delta^*(G_2); \mathbb{Z}\right)\cong
\Hred_{j-1}\left(\Delta^*(G-\{v\}) ; \mathbb{Z}\right)
\end{equation}
for all $j-1 \leq \delta-3$, i.e., for all $j\leq \delta-2$.
In particular (\ref{eqn2}) holds
for all $j\leq 2$, if $\delta=4$, and
for all $j\leq 3$, if $\delta>4$.

Assume now that $\delta=4$.
The only possible $3$-dimensional face $f\in\Delta^*(G-\{v\})$
which is not in $\Delta^*(G_1) \cap \Delta^*(G_2)$ is $f=\{v_1, v_2, v_3, v_4\}$;
this $f$ will indeed be a face of $\Delta^*(G-\{v\})$ if and only if there exist $u,w\in \V(G)-\{v,v_1, v_2, v_3, v_4\}$
so that $(u,w)\in \E(G)$.

If we can find yet another vertex $x\in \V(G)-\{v, v_1, v_2, v_3, v_4, u, w\}$ then
$g=\{x,v_1, v_2, v_3, v_4\}$ is also a face in $\Delta^*(G-\{v\})$.
Let
$$\mathcal{C}: \quad \dots \xrightarrow[]{\partial_{i+2}} C_{i+1} \xrightarrow[]{\partial_{i+1}} C_{i} \xrightarrow[]{\partial_{i}} \dots$$
be the chain complex associated with $\Delta^*(G-\{v\})$.
By considering $\partial_{3}\partial_{4}(g)$ we see that we can write $\partial_{3}(f)\in C_2$ as a $\mathbb{Z}$-linear combination of
$$ \partial_{3}(x,v_2,v_3,v_4), \partial_{3}(v_1,x,v_3,v_4), \partial_{3}(v_1, v_2, x, v_4), \partial_{3}(v_1, v_2, v_3, x) \in C_2 $$
and, since $\{x,v_2,v_3,v_4\}, \{v_1,x,v_3,v_4\}, \{v_1, v_2, x, v_4\}, \{v_1, v_2, v_3, x\}\in \Delta^*(G_1) \cap \Delta^*(G_2)$,
we deduce that (\ref{eqn2}) holds for all $j\leq 3$.

We are now left with the case where $G$ has only 7 vertices, namely, $v, v_1, v_2, v_3, v_4, u$ and $w$;
here the degree of $v$ is $7-3=4$ and the Lemma holds in this case by Lemma \ref{DegreeAtLeastn-4}(b),
so we assume now that we are not in this case.

We have just shown that (\ref{eqn2}) holds
for all $j\leq 3$ and if we combine this with
equation (\ref{eqn1}) we obtain
$$\Hred_{j}\left(\Delta^*(G_1) \cup \Delta^*(G_2); \mathbb{Z}\right)\cong
\Hred_{j-1}\left(\Delta^*(G-\{v\}) ; \mathbb{Z}\right)$$
for $j\leq 3$.

\end{proof}

\begin{thm}\label{thm2.4}
For any graph $G$, the $i$th Betti number of $\Delta(G)$ does not depend on the characteristic of $K$ for
all $0\leq i\leq 5$.
\end{thm}

\begin{proof}
Pick a counter-example $G$ with smallest number of vertices $n$.

Assume first that the degrees of the vertices of $G$ are at most $3$.
Hochster's formula implies that we need to show that
$\dim_K \Hred_{n-i-1}(\Delta(G); Z)$ is torsion-free
for all $i\leq 5$, and this is guaranteed by Lemma \ref{DegreeAtMost3}.

Assume now that there exists a vertex in $G$ with degree $\delta\geq 4$.
The Eagon-Reiner formula implies that we need to show that
$\Hred_{i-2}\left(\Delta^*(G) ; \mathbb{Z}\right)$
is torsion free for all $i\leq 5$, and this is guaranteed by Lemma \ref{lem2.4}.
\end{proof}

\section{A minimal graph ideal with characteristic-dependent Betti numbers.}\label{section4}

In this section we construct an example of a small graph ideal whose
8th Betti number differs in characteristics $0$ and $2$.
We start by recalling a well known example due to Gerald A.~Reisner.

Consider the following triangulation $\Delta^\prime$ of the real projective plane
\nobreak
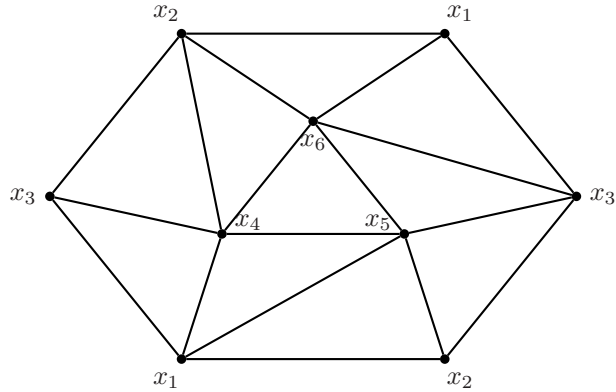
\begin{figure}[h]\label{fig1}
\begin{center}
\psset{xunit=0.7}
\psset{yunit=0.5}
\begin{pspicture}(0,0)(10,10)
\dotnode(10.000000,5.000000){v3a} \uput[0](10.000000,5.000000){$x_3$}
\dotnode(7.500000,9.330127){v1a} \uput[60](7.500000,9.330127){$x_1$}
\dotnode(2.500000,9.330127){v2a} \uput[120](2.500000,9.330127){$x_2$}
\dotnode(0.000000,5.000000){v3b} \uput[180](0.000000,5.000000){$x_3$}
\dotnode(2.500000,.669873){v1b} \uput[240](2.500000,.669873){$x_1$}
\dotnode(7.500000,.669873){v2b} \uput[300](7.500000,.669873){$x_2$}
\dotnode(6.732051,4.000000){v5} \uput[150](6.732051,4.000000){$x_5$}
\dotnode(5.000000,7.000000){v6} \uput[270](5.000000,7.000000){$x_6$}
\dotnode(3.267949,4.000000){v4} \uput[30](3.267949,4.000000){$x_4$}
\ncline{-}{v3a}{v1a}
\ncline{-}{v1a}{v2a}
\ncline{-}{v2a}{v3b}
\ncline{-}{v3b}{v1b}
\ncline{-}{v1b}{v2b}
\ncline{-}{v2b}{v3a}
\ncline{-}{v4}{v5}
\ncline{-}{v5}{v6}
\ncline{-}{v6}{v4}
\ncline{-}{v3a}{v5} \ncline{-}{v3a}{v6}
\ncline{-}{v1a}{v6}
\ncline{-}{v2a}{v6} \ncline{-}{v2a}{v4}
\ncline{-}{v3b}{v4}
\ncline{-}{v1b}{v4} \ncline{-}{v1b}{v5}
\ncline{-}{v2b}{v5}
\end{pspicture}
\end{center}
\caption{A six point triangulation of the real projective plane.}
\end{figure}

One can show that the Stanley-Reisner ring $K[\Delta]$ is Cohen-Macaulay if and only if
the characteristic of $K$ is not $2$ (cf. \S 5.3 in \cite{BH1})
so the projective dimension, and hence the Betti numbers, of $K[\Delta]$ differ in characteristics $0$ and $2$.
Specifically, when $K$ has characteristic $0$, $K[\Delta]$ has Betti number diagram
\nobreak
\begin{verbatim}
total: 1 10 15 6
    0: 1  .  . .
    1: .  .  . .
    2: . 10 15 6
\end{verbatim}
and when $K$ has characteristic $2$, $K[\Delta]$ has Betti number diagram
\begin{verbatim}
total: 1 10 15 7 1
    0: 1  .  . . .
    1: .  .  . . .
    2: . 10 15 6 1
    3: .  .  . 1 .
\end{verbatim}

We now introduce the following subdivision $\Delta$ of $\Delta^\prime$:
\begin{figure}[h]\label{fig2}
\begin{center}
\psset{xunit=0.7}
\psset{yunit=0.5}
\begin{pspicture}(0,0)(10,10)
\dotnode(10.000000,5.000000){v3a} \uput[0](10.000000,5.000000){$x_3$}
\dotnode(7.500000,9.330127){v1a} \uput[60](7.500000,9.330127){$x_1$}
\dotnode(2.500000,9.330127){v2a} \uput[120](2.500000,9.330127){$x_2$}
\dotnode(0.000000,5.000000){v3b} \uput[180](0.000000,5.000000){$x_3$}
\dotnode(2.500000,.669873){v1b} \uput[240](2.500000,.669873){$x_1$}
\dotnode(7.500000,.669873){v2b} \uput[300](7.500000,.669873){$x_2$}
\dotnode(6.732051,4.000000){v5} \uput[150](6.732051,4.000000){$x_5$}
\dotnode(5.000000,7.000000){v6} \uput[270](5.000000,7.000000){$x_6$}
\dotnode(3.267949,4.000000){v4} \uput[30](3.267949,4.000000){$x_4$}

\dotnode(8.750000000, 7.165063500){v12a} \uput[30](8.750000000, 7.165063500){$x_{12}$}
\dotnode(5.000000000, 9.330127000){v11a} \uput[80](5.000000000, 9.330127000){$x_{11}$}
\dotnode(1.250000000, 7.165063500){v10a} \uput[150](1.250000000, 7.165063500){$x_{10}$}
\dotnode(1.250000000, 2.834936500){v12b} \uput[210](1.250000000, 2.834936500){$x_{12}$}
\dotnode(5.000000000, .6698730000){v11b} \uput[270](5.000000000, .6698730000){$x_{11}$}
\dotnode(8.750000000, 2.834936500){v10b} \uput[330](8.750000000, 2.834936500){$x_{10}$}
\dotnode(5.866025500, 5.500000000){v8} \uput[30](5.866025500, 5.500000000){$x_8$}
\dotnode(4.133974500, 5.500000000){v7} \uput[150](4.133974500, 5.500000000){$x_7$}
\dotnode(5.000000000, 4.000000000){v9} \uput[330](5.000000000, 4.000000000){$x_9$}

\ncline{-}{v3a}{v1a}
\ncline{-}{v1a}{v2a}
\ncline{-}{v2a}{v3b}
\ncline{-}{v3b}{v1b}
\ncline{-}{v1b}{v2b}
\ncline{-}{v2b}{v3a}
\ncline{-}{v4}{v5}
\ncline{-}{v5}{v6}
\ncline{-}{v6}{v4}
\ncline{-}{v3a}{v5} \ncline{-}{v3a}{v6}
\ncline{-}{v1a}{v6}
\ncline{-}{v2a}{v6} \ncline{-}{v2a}{v4}
\ncline{-}{v3b}{v4}
\ncline{-}{v1b}{v4} \ncline{-}{v1b}{v5}
\ncline{-}{v2b}{v5}

\ncline{-}{v12a}{v6}
\ncline{-}{v11a}{v6}
\ncline{-}{v10a}{v4}
\ncline{-}{v12b}{v4}
\ncline{-}{v11b}{v5}
\ncline{-}{v10b}{v5}

\ncline{-}{v3a}{v8}
\ncline{-}{v2a}{v7}
\ncline{-}{v1b}{v9}

\ncline{-}{v7}{v8}
\ncline{-}{v8}{v9}
\ncline{-}{v9}{v7}

\end{pspicture}
\end{center}
\caption{A 12 point triangulation of the real projective plane.}
\end{figure}
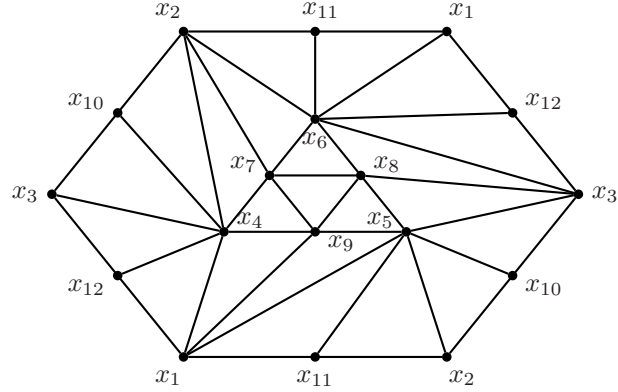

Now there exists a graph $G$ with $\Delta=\Delta(G)$, namely,
$\V(G)=\left\{ x_1, \dots, x_{12} \right\}$ and
$$\begin{array}{ll}
\E(G)=&\{
x_{1} x_{2},x_{1} x_{3},x_{1} x_{7},x_{1} x_{8},x_{1} x_{10},x_{2} x_{3}, x_{2} x_{8},x_{2} x_{9},x_{2} x_{12},x_{3} x_{7},x_{3} x_{9}, x_{3} x_{11},\\
&x_{4} x_{5},x_{4} x_{6},x_{4} x_{8},x_{4} x_{11},x_{5} x_{6},x_{5} x_{7}, x_{5} x_{12},x_{6} x_{9},x_{6} x_{10},x_{7} x_{10},x_{7} x_{11},x_{7} x_{12},\\
&x_{8} x_{10},x_{8} x_{11},x_{8} x_{12},x_{9} x_{10},x_{9} x_{11}, x_{9} x_{12}, x_{10} x_{11},x_{10} x_{12},x_{11} x_{12}\}
\end{array}
$$

The Betti numbers of $\Delta(G)$ when $K$ has characteristic $0$ are
{\nobreak
\begin{verbatim}
       total: 1 33 162 429 756 909 720 355 99 12
           0: 1  .   .   .   .   .   .   .  .  .
           1: . 33 132 228 201  93  24   3  .  .
           2: .  .  30 201 555 816 696 352 99 12
\end{verbatim}}
and when $K$ has characteristic $2$ the Betti numbers are
\nobreak
\begin{verbatim}
       total: 1 33 162 429 756 909 720 355 99 13 1
           0: 1  .   .   .   .   .   .   .  .  . .
           1: . 33 132 228 201  93  24   3  .  . .
           2: .  .  30 201 555 816 696 352 99 12 1
           3: .  .   .   .   .   .   .   .  .  1 .
\end{verbatim}
Here the 9th Betti number depends on the characteristic of $K$.

We can remove the vertex $x_2$ and some further edges to obtain a subgraph $H$ of $G$
with 11 vertices $x_1, x_3, \dots, x_{12}$ and edges
\begin{eqnarray*}
\E(G)=\{
x_{1} x_{3},x_{1} x_{7},x_{1} x_{8},x_{1} x_{10},x_{3} x_{7},x_{3} x_{9},x_{3} x_{11},x_{4} x_{5},x_{4} x_{6},x_{4} x_{8},x_{4} x_{11},x_{5} x_{6},\\
x_{5} x_{7},x_{5} x_{12},x_{6} x_{9},x_{6} x_{10},x_{7} x_{12},x_{8} x_{10},x_{8} x_{11},x_{9} x_{10},x_{9} x_{11},x_{9} x_{12},x_{11} x_{12}\}
\end{eqnarray*}

The Betti numbers of $K[\Delta(H)]$ when $K$ has characteristic $0$ are
\begin{verbatim}
       total: 1 23 103 267 442 444 259 82 11
           0: 1  .   .   .   .   .   .  .  .
           1: . 23  66  65  20   2   .  .  .
           2: .  .  37 202 422 442 259 82 11
\end{verbatim}
and when $K$ has characteristic $2$ the Betti numbers are
\begin{verbatim}
       total: 1 23 103 267 442 444 259 82 12 1
           0: 1  .   .   .   .   .   .  .  . .
           1: . 23  66  65  20   2   .  .  . .
           2: .  .  37 202 422 442 259 82 11 1
           3: .  .   .   .   .   .   .  .  1 .
\end{verbatim}
Here the 8th Betti number depends on the characteristic of $K$.

\begin{thm}\label{exampleIsMinimal}
The example $H$ above is minimal in the sense that
for any graph with at most 10 vertices,
the Betti numbers of $\Delta(G)$ do not depend on the characteristic of $G$.
\end{thm}

\begin{proof}[Computer proof:]
Pick a counter-example $G$ with minimal number of vertices $n$ for which
$\beta^K_{i,d}$ depends on $K$ for some $i,d$.
In view of Theorem \ref{thm1.0} the minimality of $G$ implies that
$d=n$ and that $\Hred_{d-i-1}(\Delta(G); \mathbb{Z})$ has torsion.
In view of Theorem \ref{thm2.4} we further assume that $i\geq 6$.

Corollary \ref{degree1} and Lemma \ref{DegreeAtLeastn-4}(b) imply that the vertices of $G$ have degree at most $n-5$ and
at least $2$. If $n< 8$, the maximal degree of vertices in $G$ is $2$ and $G$ cannot be a counter-example by
Theorem \ref{DegreeAtMost2}.

When $n=8$ we have $\beta^K_{i,d}(G)=0$ for all $i>7$ and
$\sum_{i=0}^7 (-1)^i \beta^K_{i,d}(G)= \dim_K (K[\Delta(G)])_d$
is a value of the Hilbert function of $K[\Delta(G)]$ and hence independent of $K$.
$\beta^K_{7,d}(G)$ vanishes for $d\neq 8$, Corollary \ref{cor1.2}(b) implies that $\beta^K_{7,8}(G)$
is independent of $K$ and Theorem \ref{thm2.4}
implies that $\beta^K_{i,d}(G)$ is independent of $K$ for all $0\leq i\leq 5$ and all $d$, so we conclude
that $\beta^K_{6,d}(G)$ must also be independent of $K$.

Pick any vertex $v\in\V(G)$ and let $v_1,\dots, v_s$ be its non-neighbours;
denote with $H$ the induced subgraph of $G$ with vertices $v_1,\dots, v_s$. As in the proof of
Lemma \ref{DegreeAtLeastn-4}(b), $G$ will be a minimal example
only if $\Hred_i(\Delta(H); \mathbb{Z})\neq 0$ for some $i$.
When $s=4$ Lemma \ref{lem2.0} implies that $H$ must consist of two disjoint edges and, when $s=5$,
$H$ must be one of
$$
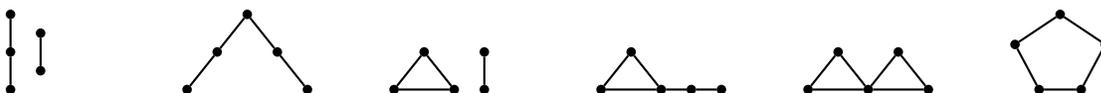
\begin{array}{l l l l l l }
\psset{xunit=0.4}
\psset{yunit=0.5}
\begin{pspicture}(0,0)(6,2)
\dotnode(1,0){v1}
\dotnode(1,1){v2}
\dotnode(1,2){v3}
\dotnode(2,0.5){v4}
\dotnode(2,1.5){v5}

\ncline{-}{v1}{v2}
\ncline{-}{v2}{v3}
\ncline{-}{v4}{v5}
\end{pspicture}
&
\psset{xunit=0.4}
\psset{yunit=0.5}
\begin{pspicture}(0,0)(6,2)
\dotnode(0,0){v1}
\dotnode(1,1){v2}
\dotnode(2,2){v3}
\dotnode(3,1){v4}
\dotnode(4,0){v5}

\ncline{-}{v1}{v2}
\ncline{-}{v2}{v3}
\ncline{-}{v3}{v4}
\ncline{-}{v4}{v5}
\end{pspicture}
&
\psset{xunit=0.4}
\psset{yunit=0.5}
\begin{pspicture}(0,0)(6,2)
\dotnode(0,0){v1}
\dotnode(1,1){v2}
\dotnode(2,0){v3}
\dotnode(3,0){v4}
\dotnode(3,1){v5}

\ncline{-}{v1}{v2}
\ncline{-}{v2}{v3}
\ncline{-}{v3}{v1}
\ncline{-}{v4}{v5}
\end{pspicture}
&
\psset{xunit=0.4}
\psset{yunit=0.5}
\begin{pspicture}(0,0)(6,2)
\dotnode(0,0){v1}
\dotnode(1,1){v2}
\dotnode(2,0){v3}
\dotnode(3,0){v4}
\dotnode(4,0){v5}

\ncline{-}{v1}{v2}
\ncline{-}{v2}{v3}
\ncline{-}{v3}{v1}
\ncline{-}{v3}{v4}
\ncline{-}{v4}{v5}
\end{pspicture}
&
\psset{xunit=0.4}
\psset{yunit=0.5}
\begin{pspicture}(0,0)(6,2)
\dotnode(0,0){v1}
\dotnode(1,1){v2}
\dotnode(2,0){v3}
\dotnode(3,1){v4}
\dotnode(4,0){v5}

\ncline{-}{v1}{v2}
\ncline{-}{v2}{v3}
\ncline{-}{v3}{v1}
\ncline{-}{v3}{v4}
\ncline{-}{v4}{v5}
\ncline{-}{v5}{v3}
\end{pspicture}
&
\psset{xunit=0.4}
\psset{yunit=0.5}
\begin{pspicture}(0,0)(6,2)
\dotnode(0.8,0){v1}
\dotnode(2.2,0){v2}
\dotnode(3,1.2){v3}
\dotnode(1.5,2){v4}
\dotnode(0,1.2){v5}

\ncline{-}{v1}{v2}
\ncline{-}{v2}{v3}
\ncline{-}{v3}{v4}
\ncline{-}{v4}{v5}
\ncline{-}{v5}{v1}
\end{pspicture}
\end{array}
$$

Assume that $n=9$.
We need to show that $\beta_{i,9}^K(G)$ is independent of $K$ or, equivalently,
by the Universal Coefficient Theorem, that $\Hred_{9-i-1}(\Delta(G); \mathbb{Z})$ has no torsion for all $5\leq i\leq 7$.
There are $5621$ unlabelled connected graphs on $9$
vertices whose degrees are $2,3,4$
\footnote{These were produced with \cite{Mc1}, see also \cite{Mc2}.}
and only $99$ of those all have vertices of degree $n-5=4$ and $n-6=3$ satisfying the conditions above.
The integral homology of all the simplicial complexes associated with these graphs was computed
\footnote{Integral homologies were computed with {\tt MOISE  - A Topology Package for Maple}
written by R.~Andrew Hicks and available from {\tt http://www.cis.upenn.edu/$\sim$rah/MOISE.html}.}
and none was found to have torsion.

Assume that $n=10$.
We need to show that $\beta_{i,10}^K(G)$ is independent of $K$ or, equivalently,
by the Universal Coefficient Theorem, that $\Hred_{10-i-1}(\Delta(G); \mathbb{Z})$
has no torsion for all $5\leq i\leq 8$.
There are $753827$ unlabelled connected graphs on $10$ vertices whose degrees are $2,3,4,5$
but (fortunately!) only $8534$ of those have all vertices of degree $n-5=5$ and $n-6=4$ satisfying the conditions above.
The integral homology of all the simplicial complexes associated with these graphs was computed and none
was found to have torsion.

\end{proof}

\begin{cor}\label{6thBetti}
For all graphs $G$, $\beta^K_{6}(\Delta(G))$ is independent of $K$.
\end{cor}
\begin{proof}
Assume we can pick a counter-example and that we pick it so that
$\beta^K_{6,j}(\Delta(G))$ depends of $K$. Lemma \ref{basicFactsOnKDependence}(d) allows us
to assume that $G$ has $j$ vertices.
Lemma \ref{BettiZeroBeyond2i} shows that, unless $7\leq j\leq 12$, $\beta^K_{6,j}(\Delta(G))=0$.
Also $\beta^K_{6,7}(\Delta(G))=\Hred_{0}(\Delta(G); K)$ is independent of $K$,
$\beta^K_{6,12}(\Delta(G))$ is independent of $K$ (by Lemma \ref{BettiZeroBeyond2i},)
$\beta^K_{6,11}(\Delta(G))$ is independent of $K$ (by Lemma \ref{Betti i,2i-1},)
and
$\beta^K_{6,8}(\Delta(G))$, $\beta^K_{6,9}(\Delta(G))$, $\beta^K_{6,10}(\Delta(G))$
are independent of $K$ (by Theorem \ref{exampleIsMinimal}.)
\end{proof}

A long search involving 2105589 graphs shows that there exist precisely \emph{four} unlabelled graphs with 11 vertices
whose Betti numbers depend on $K$, and those Betti numbers
depending on $K$ are the eighth and ninths Betti numbers (see appendix below.)

Consider now the seventh Betti number.
Assume we can pick a graph $G$  so that
$\beta^K_{7,j}(\Delta(G))$ depends of $K$ for some $j$.
Lemma \ref{basicFactsOnKDependence}(d) allows us
to assume that $G$ has $j$ vertices.
Unless $8\leq j\leq 14$, $\beta^K_{7,j}(\Delta(G))=0$.
Also $\beta^K_{7,8}(\Delta(G))=\Hred_{0}(\Delta(G); K)$ is independent of $K$,
$\beta^K_{7,14}(\Delta(G))$ is independent of $K$ (by Lemma \ref{BettiZeroBeyond2i},)
$\beta^K_{7,13}(\Delta(G))$ is independent of $K$ (by Lemma \ref{Betti i,2i-1},)
$\beta^K_{7,9}(\Delta(G))$ and $\beta^K_{7,10}(\Delta(G))$ are independent of $K$ (by Theorem \ref{exampleIsMinimal})
and
$\beta^K_{7,11}(\Delta(G))$ is independent of $K$ (by the remark above.)
Hence the only seventh Betti number which might depend on $K$ is $\beta_{7,12}^K(G)$.

\section{Appendix: graphs with 11 vertices whose Betti numbers depend on the ground field.}

There are precisely \emph{four} graphs $G_1$, $G_2$, $G_3$ and $G_4$
with 11 vertices whose Betti numbers depend on the characteristic of the ground field.
For each such graph $G_i$, $\beta_j^\mathbb{Q}(G_i)\neq \beta_j^{\mathbb{Z}/p\mathbb{Z}}(G_i)$ only when
$p=2$ and $j\in\{8,9\}$.

The edges of these graphs, together with their Betti numbers in characteristics 0 and 2 are given below.

$$\begin{array}{ll}
\E(G_1)=&\big\{ \{1, 5\}, \{1, 6\}, \{1, 8\}, \{1, 10\}, \{2, 5\}, \{2, 6\}, \{2, 9\}, \{2, 11\}, \{3, 7\},   \\
&  \{3, 8\}, \{3, 9\}, \{3, 11\},  \{4, 7\}, \{4, 8\}, \{4, 10\}, \{4, 11\}, \{5, 8\}, \{5, 9\},    \{6, 10\}, \\
& \{6, 11\}, \{7, 9\}, \{7, 10\}, \{8, 11\} \big\}
\end{array}
$$
\begin{verbatim}
total: 1 23 103 267 442 444 259 82 11
    0: 1  .   .   .   .   .   .  .  .
    1: . 23  66  65  20   2   .  .  .
    2: .  .  37 202 422 442 259 82 11

total: 1 23 103 267 442 444 259 82 12 1
    0: 1  .   .   .   .   .   .  .  . .
    1: . 23  66  65  20   2   .  .  . .
    2: .  .  37 202 422 442 259 82 11 1
    3: .  .   .   .   .   .   .  .  1 .
\end{verbatim}

$$\begin{array}{ll}
\E(G_2)=&\big\{ \{1, 4\}, \{1, 5\}, \{1, 8\}, \{1, 9\}, \{2, 5\}, \{2, 6\}, \{2, 8\}, \{2, 10\}, \{2, 11\}, \\
& \{3, 6\}, \{3, 7\}, \{3, 9\}, \{3, 10\}, \{4, 7\}, \{4, 8\}, \{4, 11\}, \{5, 9\}, \{5, 10\},  \{5, 11\}, \\
& \{6, 8\}, \{6, 9\}, \{6, 11\}, \{7, 10\}, \{7, 11\} \big\}
\end{array}
$$
\begin{verbatim}
total: 1 24 104 257 419 425 252 81 11
    0: 1  .   .   .   .   .   .  .  .
    1: . 24  73  80  30   4   .  .  .
    2: .  .  31 177 389 421 252 81 11

total: 1 24 104 257 419 425 252 81 12 1
    0: 1  .   .   .   .   .   .  .  . .
    1: . 24  73  80  30   4   .  .  . .
    2: .  .  31 177 389 421 252 81 11 1
    3: .  .   .   .   .   .   .  .  1 .
\end{verbatim}

$$\begin{array}{ll}
\E(G_3)=& \big\{\{1, 4\}, \{1, 5\}, \{1, 8\}, \{1, 9\}, \{2, 5\}, \{2, 6\}, \{2, 8\}, \{2, 10\}, \{2, 11\}, \\
&  \{3, 6\}, \{3, 7\}, \{3, 9\}, \{3, 10\}, \{4, 7\}, \{4, 8\}, \{4, 11\}, \{5, 9\}, \{5, 10\},    \{5, 11\}, \\
&  \{6, 8\}, \{6, 9\}, \{6, 11\}, \{7, 10\}, \{7, 11\}, \{9, 11\} \big\}
\end{array}
$$
\begin{verbatim}
total: 1 25 107 255 406 411 246 80 11
    0: 1  .   .   .   .   .   .  .  .
    1: . 25  80  97  46  10   1  .  .
    2: .  .  27 158 360 401 245 80 11

total: 1 25 107 255 406 411 246 80 12 1
    0: 1  .   .   .   .   .   .  .  . .
    1: . 25  80  97  46  10   1  .  . .
    2: .  .  27 158 360 401 245 80 11 1
    3: .  .   .   .   .   .   .  .  1 .
\end{verbatim}

$$\begin{array}{ll}
\E(G_4)=&\big\{ \{1, 4\}, \{1, 5\}, \{1, 7\}, \{1, 8\}, \{2, 5\}, \{2, 6\}, \{2, 8\}, \{2, 10\}, \{2, 11\},    \\
& \{3, 6\}, \{3, 7\}, \{3, 9\}, \{3, 10\}, \{4, 8\}, \{4, 9\}, \{4, 10\}, \{4, 11\}, \{5, 7\},    \{5, 9\}, \\
&  \{5, 11\}, \{6, 8\}, \{6, 9\}, \{7, 10\}, \{9, 11\}, \{10, 11\} \big\}
\end{array}
$$
\begin{verbatim}
total: 1 25 105 247 396 406 245 80 11
    0: 1  .   .   .   .   .   .  .  .
    1: . 25  80  95  40   6   .  .  .
    2: .  .  25 152 356 400 245 80 11

total: 1 25 105 247 396 406 245 80 12 1
    0: 1  .   .   .   .   .   .  .  . .
    1: . 25  80  95  40   6   .  .  . .
    2: .  .  25 152 356 400 245 80 11 1
    3: .  .   .   .   .   .   .  .  1 .
\end{verbatim}

\section*{Acknowledgment}

I would like to thank the two anonymous referees of this paper for their very useful corrections and suggestions.


\begin{thebibliography}{BH1}

\bibitem[B]{B}
Dave Bayer,
\emph{Monomial ideals and duality,} Lecture notes, Berkeley 1995-96,
available by anonymous ftp at {\tt math.columbia.edu/pub/bayer/monomials duality/monomials.ps}.

\bibitem[BS]{BS}
Dave Bayer and  Bernd Sturmfels.
\emph{Cellular resolutions of monomial modules.}
J.~Reine Angew.~Math.~502 (1998), pp.~123--140.

\bibitem[BH1]{BH1}
Winfried Bruns and J\"{u}rgen Herzog.
\emph{Cohen-Macaulay Rings},
Cambridge University Press,
Cambridge (1993)

\bibitem[BH2]{BH2}
Winfried Bruns and J\"{u}rgen Herzog.
\emph{On multigraded resolutions.}
Mathematical Proceedings of the Cambridge Philosophical Society 118 (1995), no.~2, pp.~245--257.

\bibitem[ER]{ER}
John A.~Eagon and Victor Reiner.
\emph{Resolutions of Stanley-Reisner rings and Alexander duality.}
J. Pure Appl. Algebra 130 (1998), no. 3, pp.~265--275.

\bibitem[H]{H}
Melvin Hochster.
\emph{Cohen-Macaulay rings, combinatorics, and simplicial complexes.}
Ring theory, II (Proc. Second Conf., Univ. Oklahoma, Norman, Okla., 1975),
pp.~171--223. Lecture Notes in Pure and Appl. Math., Vol. 26, Dekker, New York, 1977.

\bibitem[J]{J}
Sean Jacques.
\emph{The Betti numbers of graph ideals.}
PhD Thesis, The University of Sheffield, 2004.
arXiv.math.AC/0410107.

\bibitem[JK]{JK}
Sean Jacques and Mordechai Katzman.
\emph{The Betti numbers of forests.}
arXiv.math.AC/0501226

\bibitem[Mc1]{Mc1}
Brendon D.~McKay.
``{\tt geng.}''
A program for generating graphs available at
{\tt http://cs.anu.edu.au/$\sim$bdm/nauty/}

\bibitem[Mc2]{Mc2}
Brendon D.~McKay.
\emph{Isomorph-free exhaustive generation.}
J.~Algorithms, {\bf 26} (1998) 306--324.

\bibitem[M]{M}
William S.~Massey.
\emph{A basic course in algebraic topology.}
Graduate Texts in Mathematics, 127.
Springer-Verlag, New York, 1991.

\bibitem[S]{S}
Richard P.~Stanley.
\emph{Combinatorics and commutative algebra.}
Progress in Mathematics, {\bf 41.} Birkh\"{a}user Boston, Inc., Boston, MA, 1983.

\bibitem[T]{T}
Diana Taylor.
\emph{Ideals generated by monomials in an $R$-sequence.}
Thesis, Chicago University, 1966.


\bibitem[TH]{TH}
Naoki Terai and Takayuki Hibi.
\emph{Some results on Betti numbers of Stanley-Reisner rings.}
Proceedings of the 6th Conference on Formal Power Series and Algebraic Combinatorics (New Brunswick, NJ, 1994).
Discrete Math. {\bf 157} (1996), no.~1-3, 311--320.

\end{thebibliography}
\end{document}